\documentclass[a4paper]{article}
\usepackage{latexsym}
\usepackage{setspace}
\usepackage{amssymb,amsmath}
\usepackage{amsthm}
\usepackage[english]{babel}
\usepackage[applemac]{inputenc} 
\usepackage{epsfig}
\usepackage{graphicx}

\numberwithin{equation}{section}
\newtheorem{theorem}{Theorem}[section]
\newtheorem{definition}[theorem]{Definition}

\newtheorem{proposition}[theorem]{Proposition}

\newtheorem{remark}[theorem]{Remark}
\newcommand{\cvd}{$\hfill \sqcap \hskip-6.5pt \sqcup$} %(quadratino bianco)
\newcommand{\INT}{\int_{\Om}}
\newcommand{\acca}{H^1_0(\Omega)}
\newcommand{\Om}{\Omega}
\newcommand{\erre}{\mathbb{R}}

\newcommand{\integrale}[1]{\int _{\Omega} #1 \,dx}

\newcommand{\p}{$\mathcal{P}$}

\newcommand{\sublim}[2]{\substack{#1 \\ #2}}

\renewcommand{\rho}{\varrho}
\renewcommand{\theta}{\vartheta}
\newcommand{\be}{\begin{equation}}
\newcommand{\ee}{\end{equation}}

\newcommand{\na}{\nabla}
\newcommand{\ba}{\begin{array}}
\newcommand{\ea}{\end{array}}

\begin{document}
\title{Minimax solutions for a problem with sign changing nonlinearity and
lack of strict convexity}
\date{ }
\maketitle
\centerline {\large {Paola Magrone}}
\centerline{\large {Dipartimento di Architettura,
Università degli Studi Roma Tre}}
\centerline{\large{ 
Via della Madonna dei Monti 40, Roma, Italia}}

%00147 Roma, Italia}
%\email{magrone@mat.uniroma3.it}
\begin{abstract}
A result of existence of a nonnegative and a nontrivial solution is proved via critical point theorems for non smooth 
functionals. The equation considered presents a convex part and
a nonlinearity which changes sign.\footnote{ AMS Subject Classification 35J65, 58E05; keywords: non strict convexity, sign changing, Linking theorem}
\end{abstract}
 %\subjclass{35J65, 58E05}
%\keywords{non strict convexity, sign changing, Linking theorem, nontrivial solution}

%\noindent {\bf AMS Subject Classification}: 49J40, 35J85, 58E05.

\section{Introduction and main results}\label{sect:intro}
Let us consider the problem
\be
\tag{\p}
\left\{
\begin{array}{ll}
- \mathrm{div}(\Psi'(\na u))=
\lambda u+b(x)|u|^{p-2}u &\mbox{in $\Om$,}\\
\noalign{\medskip}
u=0 &\mbox{on $\partial\Om$,}
\end{array}\right.
\ee where $\lambda$ is a real parameter, $\Om$ is a bounded open
subset of $\erre^N,$ $N\geq 2,$ $b(x)\in\overline C(\Om)$ changes sign in $\Om.$
Finally $2<p<2^*=\frac{2N}{N-2},$ and we will assume that $\Psi:\erre^N\to\erre$ is a
convex function of class~$C^1$ satisfying the following
conditions:
\vskip0,2cm
$\begin{array}{ll}
(\Psi_1) &
\displaystyle{
\lim_{\xi\to 0}\frac{\Psi(\xi)}{|\xi|^2}=\frac12;} \\
& \\
(\Psi_2) &
\displaystyle{ \exists \mu>0 :\ \mu |\xi|^2\leq \Psi(\xi)\leq \frac1\mu |\xi|^2\quad \mbox{for every $\xi\in\mathbb{R}^N$}
;} \\
& \\
(\Psi_3) &
\displaystyle{ \lim_{|\xi|\to \infty} \frac{ \Psi'(\xi)\cdot\xi-2\Psi (\xi)}{|\xi|^2}=0
;}

\end{array}$
%(\Psi_3) &
%\displaystyle{\Psi(\xi)\leq \frac12|\xi|^2
%\qquad\mbox{for every $\xi\in\mathbb{R}^N$}.}
%\end{array}$
%

\vskip0,2cm
\noindent
Moreover the function $b(x)$ has to be strictly positive in a non zero measure set, and the zero set  must be "thin", in other words
$b(x)$ must satisfy the following conditions:
\vskip0,2cm
\noindent
$\begin{array}{ll}
(b_1) &
\displaystyle{\Omega^+ :=\{x\in\Omega \ :\ b(x)>0\}\ \hbox{ is a nonempty open set}} \\
& \\
(b_2) &
\displaystyle{ \Omega^0 :=\{x\in\Omega \ :\ b(x)=0\}\ \hbox{ has zero measure }}
\end{array}$
\vskip0,2cm
\noindent

Conditions $(b_1)$ and $(b_2)$ imply that
 $ b^+(x)=b(x)+b^-(x)\not\equiv 0$ and that, since $b$ is continuous, the set $\Om^0$ is closed in $\Om$.
\vskip0,3cm
Let us also denote by $(\lambda_k)$ the eigenvalues of $-\Delta$ with homogeneous Dirichlet boundary
condition.\vskip0,1cm\noindent
In the model case $\Psi(\xi)=\frac12|\xi|^2$, there is a wide
literature on problem (\p).\\
To cite only some of the existing results, in \cite{AT1} the authors found positive solutions to (\p) in case that
$\lambda_1<\lambda<\Lambda^*,$ with $\Lambda^*$ suitably near to $\lambda_1.$ In the following many other papers (\cite{AdP}, \cite{AT1}, \cite{AT2}, \cite{B}, \cite{BaNa}) were devoted to prove existence of (possibly infinitely many) solutions for 
$\lambda\in[\lambda_1,\Lambda^*]$ or also for every $\lambda,$ in case the nonlinearity satisfies some oddness assumption. A result concerning all $\lambda$ different from the eigenvalues of the Laplacian under some quite general assumptions can be found in  \cite{RTT}, while in \cite{gmm} the authors proved a result of existence of a nontrivial solution (possibly changing sign) for every $\lambda.$
\vskip0,2cm
On the other hand, only a small literature is available when dealing
with equations with a non strictly convex principal part. In this framework, in \cite{Deg1} the author applies non smooth variational methods in presence of subcritical, positive, nonlinearities; while using similar techniques a nonlinearity with critical growth was considered in \cite{m}. 
\vskip0,2cm\noindent
The aim of this paper is to extend to the setting of non strictly convex functionals some of the results contained in \cite{AT1} (existence of a positive solution for $\lambda<\lambda_1$) and \cite{gmm} 
(existence of a nontrivial solution for any $\lambda.$)\vskip0,2cm\noindent
Problem (\p) can be treated by variational techniques.
Indeed, weak solutions $u$ of (\p) can be found as critical
points of the~$C^1$~functional $J:\acca\to\erre$ defined as
\be
\label{J}
J(u)=\int_{\Om}\Psi(\na u)\,dx -\frac\lambda{2}\int_{\Om}u^2\,dx -
\frac1{p}\int_{\Om}b(x)|u|^{p}\,dx .
\ee
\vskip0,2cm\noindent
The key point here is that, although $\Psi$ shares some properties
with this typical case, there is no assumption of strict convexity
with respect to $\xi$.
\par
For instance, one could consider 
\begin{equation}
\label{eq:ex}
\Psi(\xi) = \psi(\xi_1) + \frac12 \sum_{j=2}^N \xi_j^2\,,
\end{equation}
where
\[
\psi(t) = \left\{
\begin{array}{ll}
\frac12t^2 & \mbox{if $|t|<1$,}\\
\noalign{\medskip}
|t|-\frac12 & \mbox{if $1\leq |t|\leq 2$},\\
\noalign{\medskip}
\frac12|t|^2 - |t| + \frac32 & \mbox{if $|t|>2$.}
\end{array}\right.
\]

If we look at the principal part of $J$ as the energy stored in
the deformation $u$, this means that the material has a plastic
behavior when $1\leq |D_1 u|\leq 2$.
We refer the reader to \cite[Chapter 6]{wu} for a discussion
of several models of plasticity.
\par

As shown in \cite{Deg1,m}, it may happen that Palais Smale sequences, even if bounded in $\acca$-norm, do not admit any subsequence which converges strongly in
this norm. And there is no way to prevent the interaction between the area where $\Psi$ loses strict convexity and
the values of $\nabla u$. A possibile strategy is to look for compactness in a weaker norm ($L^{2^*}$).
\vskip0,1cm\noindent
Let us introduce the following notations: let $k\geq 1$ be such that $\lambda_k \leq \lambda < \lambda_{k+1}$
and let $e_1,\ldots,e_k$ be eigenfunctions of $-\Delta$
associated to $\lambda_1,\ldots,\lambda_k$, respectively.
Finally, let $E_-=\mathrm{span}\{e_1,...,e_k\}$ and
$E_+=E_-^\perp$.
The main result of this paper are the following:
\begin{theorem}
\label{esistenza} Let $N\geq2$ and let $ \Psi:\erre^N\to\erre$ be
a convex function of class $C^1$ satisfying $(\Psi_1),(\Psi_2)$.
Moreover let the function $b(x)$ verify $(b_1).$
Then, for every
$\lambda\in]0,\lambda_1[$, problem \em(\em\p\em)\em admits a
nontrivial and nonnegative weak solution $u\in\acca$.
%\be\label{lambdapsi}
%(1-\delta)|\xi|^2\leq\Psi'(\xi)\xi\leq2\Psi(\xi)\leq(1+\delta)|\xi|^2
%\ee
\end{theorem}
\begin{theorem}\label{teocambiosegno}
Let $N\geq2$ and let $ \Psi:\erre^N\to\erre$ be a convex function of
class $C^1$ satisfying $(\Psi_1),\ (\Psi_2)$ and let
$\lambda\geq\lambda_1$. Moreover let the function $b(x)$ verify $(b_1),$ and the following assumptions:
\be \label{ipotesib-seconda}
\int_\Om b(x)|v|^p\geq0 \qquad\forall v\in E_-.
\ee
\be \label{ipotesib-terza}
\exists e\in E_-^\perp\setminus\{0\}\ :\
\int_\Om b(x)|v|^p\,dx\geq C \int_\Om |v|^p \,dx\qquad\forall v\in E_-\oplus span \{e\}.
\ee
Then problem \em(\em\p\em)\em\
admits a nontrivial weak solution $u\in\acca$.
%\be\label{lambdapsi}
%(1-\delta)|\xi|^2\leq\Psi'(\xi)\xi\leq2\Psi(\xi)\leq(1+\delta)|\xi|^2
%\ee
\end{theorem}

%%%%%%%%%%%%%%%%%%%%%%%%%%%%%%%%%%%%%%%%%%%%%%%%%%%%%%%%%%%%%%%%%%%%%
%%%%%%%%%%%%%%%%%%%%%%%%%%%%%%%%%%%%%%%%%%%%%%%%%%%%%%%%%%%%%%%%%%%%%
%%%%%%%%%%%%%%%%%%%%%%%%%%%%%%%%%%%%%%%%%%%%%%%%%%%%%%%%%%%%%%%%%%%%%
\begin{remark}
\label{proprieta-convesse}
Arguing as in section 2 of \cite{m} we can deduce the following properties for $\Psi,$ up to modifying the constant $\mu:$
\begin{align}
\label{psi'xi}
\Psi'(\xi)\cdot\xi&\geq \mu |\xi|^2
\qquad\forall \xi\in\erre^N\,, \\
\label{eqpsi2+psi4,1}
|\Psi'(\xi)|&\geq \mu |\xi|^{\phantom{2}}
\qquad\forall \xi\in\erre^N
\end{align}

\be\label{Psi'}
|\Psi'(\xi)|\leq \frac1\mu |\xi|\qquad\forall \xi\in\erre^N
\ee
Furthermore $(\Psi_3)$ yields that  $\forall \sigma>0,$ $\exists M_\sigma\in \erre:$

\be\label{psi'-psi}
 \Psi'(\xi)\xi-2\Psi(\xi)\leq  \sigma |\xi|^2+M_\sigma
\ee
\end{remark}

%%%%%%%%%%%%%%%%%%%%%%%%%%%%%%%%%%%%%%%%%%%%%%%%%%%%%%%%%%%%%%%%%%%%
%%%%%%%%%%% variational framework %%%%%%%%%%%%%%%%%%%%%%%%%%%%%%%%%
%%%%%%%%%%%%%%%%%%%%%%%%%%%%%%%%%%%%%%%%%%%%%%%%%%%%%%%%%%%%%%%%%%%%
\section{The variational framework}
Let $\Omega$ be a bounded open subset of $\erre^N$, $N\geq 2$,
with Lipschitz boundary and let $\lambda\in\erre$. Let us define the following
functional $J:\acca\to \erre$
$$
J(u)= \int_{\Om}\Psi(\na u)\,dx -\frac\lambda{2}\int_{\Om}u^2\,dx -
\frac1{p}\int_{\Om}b(x)|u|^{p}\,dx.
$$
By $(\Psi_1)$, $(\Psi_2)$ the functional $J$ is of class $C^1$ on $\acca.$ We wish to apply variational methods to functional $J,$ but,
as already mentioned, it is well known
that the Palais Smale (PS) condition for a functional which is not strictly convex is not satisfied on $\acca.$
So it is convenient to extend the functional $J$ to $L^{2^*}$ with value $+\infty$ outside
$\acca.$
\vskip0,1cm\noindent
In other words we define the convex, lower semicontinuous functional (still denoted $J$)
\vskip0,1cm\noindent
$$J:L^{2^*}(\Omega)\longrightarrow [0,+\infty]$$
\vskip0,2cm\noindent
\be\label{funzionale-esteso}
J(u) = \left\{
\begin{array}{l}
\displaystyle{
 \int_{\Om}\Psi(\na u)\,dx -\frac\lambda{2}\int_{\Om}u^2\,dx -
\frac1{p}\int_{\Om}b(x)|u|^{p}\,dx}
\,\text{ if $u\in \acca$}\,,\\
\noalign{\medskip}
+\infty
\qquad\text{if $u\in L^{2^*}(\Omega)\setminus \acca$}\,
\end{array}
\right.
\ee
\vskip0,2cm\noindent
This setting will allow us to recover PS condition.
\vskip0,1cm\noindent
This functional can be written as $J=J_0+J_1,$ where
$$
J_0=\int_{\Om}\Psi(\na u)\,dx,
$$
is proper, convex and l.s.c., while
$$
J_1=-\frac\lambda{2}\int_{\Om}u^2\,dx -\frac1{p}\int_{\Om}b(x)|u|^{p}\,dx,
$$
is of class $C^1.$ We will use the following definitions (\cite{sz}, \cite{Deg1}) of critical point and PS sequence for functionals
 of the type $J=J_0+J_1$:
\begin{definition}\label{punto-critico}
Let $X$ be a real Banach space, $u\in X$ is a critical point for $J$ if $J(u)\in\erre$ and $-J_1'(u)\in\partial J_0,$ where $\partial J_0$ is
 the subdifferential of $J_0$ at $u.$
\end{definition}
\vskip0,3cm\noindent
\begin{definition}\label{pscondition}
Let $X$ be a real Banach space and let $c\in\erre$. We say that $u_k$ is a Palais Smale sequence at level $c$ ($(PS)_c$ sequence for short) for $J$
if $J(u_k)\to 0$ and there exists $\alpha_k\in \partial J_0$ with
$(\alpha_k+J_1'(u_k))\to 0$ in $ X^*.$
\end{definition}
\vskip0,2cm\noindent

The following proposition (see \cite{Deg1}) assures that the critical points of the extendend functional already defined gives the solutions of our problem.
\begin{proposition}\label{deg-bergamo1}
Let $u\in L^{2^*}(\Omega,\erre^N).$ Then $u$ is a critical point of $J$ if and only if $u\in \acca$ and $u$ is a weak solution of (\p).
\end{proposition}
\vskip0,1cm\noindent
\textbf{Proof}
Let $v\in L^{2^*}.$ Then $v\in \partial J_0,$ if and only if $u\in\acca$ and
$$
-div (\Psi'(\na u))=v
$$
that is a reformulation of definition \ref{punto-critico}.
\par\medskip\noindent
\cvd
\par\medskip\noindent

Moreover we will apply the compactness result contained in \cite{Deg1}, which we recall.
\vskip0,1cm\noindent
Let us define the functional  $\mathcal{E}:W^{1,2}_0(\Om,\erre^N)\to \erre$ as
$$
\mathcal{E}(u)=\INT \Psi(\na u)\,dx
$$
\begin{theorem}
Assume that $\Om$ is bounded. If $\{u_h\}$ is weakly convergent to $u$ in $W^{1,2}_0(\Om,\erre^N)$ with $\mathcal{E}(\{u_h\})\to \mathcal{E}(\{u\}),$ then $u$ is strongly
convergent to $u$ in $L^{2^*}(\Om).$
\end{theorem}

%%%%%%%%%%%%%%%%%%%%%%%%%%%%%%%%%%%%%%%%%%%%%%%%%%%%%%%%%%%%%%%%%%%%%%%%%%%%
%%%%%%%%%%%%%%%%%%%%%  main results   %%%%%%%%%%%%%%%%%%%%%%%%%%%%%%%%%%%%%%
%%%%%%%%%%%%%%%%%%%%%%%%%%%%%%%%%%%%%%%%%%%%%%%%%%%%%%%%%%%%%%%%%%%%%%%%%%%%
\section{Proof of main results}
Since $\Psi'(0)=0$, of course $0$ is a solution of (\p). Therefore we are interested in {\em nontrivial} solutions.
In order to find nonnegative solutions of (\p),we consider the modified functional
$\overline{J}:L^{2^*}(\Om)\to ]-\infty,\ +\infty]$ defined as
\[
\overline{J}(u)= \left\{
\begin{array}{l}
\displaystyle{
\int_{\Om}\Psi(\na u)\,dx -\frac\lambda{2}\int_{\Om}(u^+)^2\,dx -
\frac1{p}\int_{\Om}b(x)(u^+)^{p}\,dx}
\text{   if $u\in \acca$}\,,\\
\noalign{\medskip}
+\infty
\qquad\text{if $u\in L^{2^*}(\Omega)\setminus \acca$}\,
\end{array}
\right.
\]
Of course, $\overline{J}$ is also convex and lower semicontinuous.
\begin{proposition}
\label{natmod}
Let $\Psi:\erre^N\to\erre$ be a convex function of class $C^1$
satisfying $(\Psi_2)$ with $\mu>0$, and (\ref{eqpsi2+psi4,1}).
Then each critical point $u\in L^{2^*}$ of $\overline{J}$
is a nonnegative solution of {\em (\p)}.
\end{proposition}
\textbf{Proof}
Since by Proposition \ref{deg-bergamo1} we already know that the critical points of $J$ are solutions of our problem,
 it is only left to prove that the modified functional will give nonnegative solutions. By $(\Psi_2)$
one has
\[
\mu \integrale{|\na u^-|^2}\,dx  \leq
\integrale{\Psi'(\na u)\cdot(-\na u^-)\,dx} =\]
\[
=\lambda\,\integrale{u^+(-u^-)\,dx} + \integrale{(u^+)^{p-1}(-u^-)\,dx} =
0 \,,
\]
\vskip0,1cm\noindent
whence the assertion.

\cvd
\par\medskip\noindent

\begin{remark}
From now on, to simplify notations, we will keep on using the functional $J$ instead of $\overline J,$ since it is
understood what has been proved in Proposition \ref{natmod}.
\end{remark}

\textbf{Proof of Theorem \ref{esistenza}}
\par\medskip\noindent
We aim to apply to $J$ a nonsmooth version of Mountain Pass Theorem \cite{sz}.
First of all, let us observe that, by $(\Psi_1)$, we have
\[
\frac{\int_{\Om}\Psi(\na u)\,dx}{
\int_{\Om} |\na u|^2\,dx} \to \frac{1}{2}
\qquad\text{as $u\to 0$ in $L^{2^*}$.}
\]
Then, as in the case $\Psi(\xi)=\frac{1}{2}|\xi|^2$
treated in \cite{AT1,gmm}, we deduce that there exist
$\rho>0$ and $\alpha>0$ such that $J(u)\geq\alpha$
whenever $\|u\|=\rho$.
On the other hand, there exists $e \in L^{2^*}$ with $e\geq 0$
a.e. in $\Om$ such that
\begin{gather*}
\lim_{t\to+\infty}  J(te) = -\infty\,, \\
\end{gather*}
again, this is proved in \cite{AT1} in the case
$\Psi(\xi)=\frac{1}{2}|\xi|^2$, but by $(\Psi_2)$ the
assertion is true also in our case.

\vskip0,2cm\noindent
By the Mountain Pass theorem, there exist a sequence $(u_k)$ in
$L^{2^*}$ and a sequence $(w_k)$ in $L^{(2^*)'}(\Omega)$ strongly
convergent to $0$ such that (see definition \ref{pscondition})
\begin{multline}\label{ps-diseq1}
\INT \Psi'(\na u_k)(\na v-\na u_k)\,dx\geq \lambda \INT u_k(v-u_k)\,dx+\INT b(x)|u_k|^{p-1}(v-u_k)\,dx\\
+\INT w_k(v-u_k)\,dx
\quad\forall v\in L^{(2^*)'}
\end{multline}
Taking $v=0$ and $v=2u_k$ as tests in the previous inequality yield
\be\label{ps-diseq2}
\INT \Psi'(\na u_k)\na u_k\,dx = \lambda \INT (u_k)^2\,dx+\INT b(x)|u_k|^{p}\,dx+\INT w_k u_k\,dx
\quad\forall v\in L^{(2^*)'} .
\ee
Furthermore also the following relation holds:
\be\label{ps2}
\lim_{k\to\infty} \left(
\INT \Psi(\na u_k)\,dx
- \frac{\lambda}{2}\int_{\Om} (u_k)^2\,dx
- \frac{1}{p}\int_{\Om}b(x)(u_k)^{p}\,dx\right) = c >\alpha.
\ee

%\be\label{ps1}
%\int_{\Om}\Psi'(\na u_k)\cdot\na v\,dx
%- \lambda\int_{\Om}u_k v\,dx
%- \int_{\Om}b(x)(u_k)^{p-1} v\,dx = \INT w_k v
%\qquad\forall v\in L^{2^*}\,,
%\ee

Let us write the expression $p J(u_k)- J'(u_k)u_k,$ which is boundend by assumptions (\ref{ps-diseq2}), (\ref{ps2}):
\begin{multline}\label{pJ-J}
p\int_{\Om}\Psi(\na u_k)\,dx-\frac{p}{2}\lambda\int_{\Om}(u_k)^2\,dx
- \int_{\Om}b(x)(u_k)^{p}\,dx-\int_{\Om}\Psi'(\na u_k)\cdot\na u_k\,dx\\
+\lambda\int_{\Om}(u_k)^2\,dx + \int_{\Om}b(x)(u_k)^{p}\,dx=\\
\int_{\Om}(p-2)\Psi(\na u_k)\,dx+\INT\left[2\Psi(\na u_k) -
\Psi'(\na u_k)\cdot\na u_k\right]\,dx-\lambda\left(\frac{p}{2}-1\right)\INT(u_k)^2\,dx=\\
(p-2)c-\INT  w_k u_k\,dx+C
\end{multline}
By (\ref{psi'-psi}) and $(\Psi_2)$ one gets

\be\label{stima-norme}
\mu(p-2-\sigma) \INT |\na u_k|^2\,dx -\lambda\left(\frac{p}{2}-1\right)\lambda\int_{\Om}(u_k)^2\,dx \leq
pc-\INT w_k u_k +C
\ee
so
\be\label{stima-norme1}
\mu(p-2-\sigma) \INT |\na u_k|^2\,dx\leq
\lambda\left(\frac{p}{2}-1\right)\INT (u_k)^2\,dx+C\,
\ee
where the quantity $(p-2-\sigma)$ is strictly positive since $\sigma$ is arbitrarily small.
Our aim is to prove the boundedness of the $H^1_0$ norm of the Palais Smale sequences, so arguing by contradiction, let us assume that
$$
||u_k||\to\infty\quad \hbox{as }k\to+\infty.
$$
Dividing (\ref{ps2}) by $||u_k||^p$ yields
$$
\liminf\left\{\frac{p\int_{\Om}\Psi(\na u_k)}{||u_k||^p}\,dx-\frac{\lambda p}{2}\frac{\int_{\Om}(u_k)^2\,dx}{||u_k||^p}\,dx
- \frac1p\int_{\Om}b(x)\left(\frac{u_k}{||u_k||}\right)^{p}\,dx\right\}=0.
$$
Since $p>2$ and $(\Psi_2)$ holds, the first two terms go to zero. So
\be\label{limsup_fa_0}
\limsup\left(\int_{\Om}b(x)\left(\frac{u_k}{||u_k||}\right)^{p}\,dx\right)=0.
\ee
Since $b$ is bounded, by Lebesgue dominated convergence Theorem we can take the limit and
deduce that
\be\label{lim_b(x)_fa_0}
\lim_k b(x)\left(\frac{u_k}{||u_k||}\right)^p=0.
\ee
This yields that 
$$
\left(\frac{u_k}{||u_k||}\right)\to u_0
$$
strongly in $L^p$ and weakly in $\acca.$ Arguing by contradiction let us suppose that $u_0\equiv 0.$
Dividing (\ref{stima-norme1}) by $||u_k||^2$ yields
\be\label{stimanorme2}
\mu(p-2-2\sigma) \leq
\lambda\left(\frac{p}{2}-1\right)\frac1{||u_k||^2}\INT (u_k)^2\,dx +\frac{C}{||u_k||^2}
\ee
the right hand side goes to zero,  which leads to a contradiction since $p-2-2\sigma>0$ and $\mu>0,$ so $u_0$ must
not be identically zero.
\vskip0,2cm\noindent
%Let us observe that the following inequality holds (POTREBBE SERVIRE)
%$$
%\left(\frac{p-2}{2p}\right)\INT b(x)(u_k)^p< J(u_k) -\frac12 \INT w_k u_k.
%$$

Now let  $\phi \in C_0^{\infty}(\Omega^+)$ be a compact support function, $\phi\geq 0$ and $\phi\not\equiv 0.$
 Let us use the function $t\phi v,$ $v\in\acca$ as a test in \ref{ps-diseq1}.
 $$
 \int_{\Omega^+} \Psi'(\na u_k) (t\phi\na v+tv\na\phi-\na u_k)\geq
 \lambda \int_{\Omega^+} u_k(tv\phi-u_k)+ \int_{\Omega^+} b(x) (u_k)^{p-1}(tv\phi-u_k)
 $$
 $$ \int_{\Omega^+} w_k(tv\phi-u_k)\quad\forall v\in\acca $$
 
 Then let us divide the previous inequality by $t$ and then let $t$ go to $+\infty:$
 
\begin{multline}\label{disuguaglianza_b+}
 \int_{\Omega^+} \Psi'(\na u_k) (\phi\na v) +\Psi' (\na u_k)v\na\phi\geq
 +\lambda\int_{\Omega^+} (u_k)^2v\phi\ +\int_{\Omega^+} b^+(x) (u_k)^{p-1}v\phi\ +\\
 \\
+ \INT w_k v\phi \quad\qquad\forall v\in\acca 
 \end{multline}

On the other hand, if $t\to-\infty,$ one gets the opposite inequality, so we can deduce
that the equality holds in the last expression, that is
\vskip0,1cm\noindent 

\begin{multline}
\label{uguaglianza_b+}
 \int_{\Omega^+} \Psi'(\na u_k) (\phi\na v) +\Psi' (\na u_k)v\na\phi =
\\
 +\lambda\int_{\Omega^+} (u_k)^2v\phi+\int_{\Omega^+} b^+(x) (u_k)^{p-1}v\phi+\INT w_k v\phi \quad\forall v\in\acca .
 \end{multline}

\vskip0,3cm\noindent 

Now let us choose $v=u_k$ and divide both handsides of  (\ref{uguaglianza_b+}) by $||u_k||^p.$  It is easily seen that the terms containing $\lambda$ and $w_k$ go to $0$ as $k\to+\infty.$ Then

$$
\displaystyle{\int_{\Omega^+}\frac{\Psi'(\na u_k)\na u_k\phi}{||u_k||^p}}
$$ 
goes to $0$ since $p>2$ and (\ref{Psi'}) holds. 
\vskip0,2cm\noindent
On the other hand, by (\ref{Psi'}), since $p>2$ and $\phi$ is of class $C^\infty$ in $\Omega^+$ bounded,

$$
\displaystyle{\frac1{||u_k||^p}\int_{\Omega^+}\Psi'(\na u_k)u_k\na\phi}\leq
C\displaystyle{\frac{||u_k||}{||u_k||^{p-1}} \frac{||u_k||_{L^2}}{||u_k||} }
$$ 

The term  $ \frac{||u_k||_{L^2}}{||u_k||}$ is bounded, while $\frac{||u_k||}{||u_k||^{p-1}} $
converges to $0$.
\vskip0,3cm\noindent
By (\ref{uguaglianza_b+})  We can conclude that  
$$
\int_{\Omega^+} \frac1{||u_k||^p}b^+(x) (u_k)^{p}\phi\ \mapsto 0 \hbox{ as  } k\to\infty.
$$
Applying Fatou's Lemma yields 
$$
\liminf \int_{\Omega^+}\frac1{||u_k||^p} b^+(x) (u_k)^{p}\phi \leq 0
$$
and since the integrand is nonnegative, this means that $\frac{u_k^p}{||u_k||^p}$ must tend to $0$ in $\Om^+$
as $k\to\infty,$ but this
is in contradiction with the fact that it was already proved that it 
 converges to a nonzero function $u_0.$ 
 \vskip0,1cm\noindent
 Arguing in the same way, choosing now a compact support function $\eta\in C^\infty_0(\Om^-),$ yields that 
 $\frac{u_k^p}{||u_k||^p}\to 0$  as $k\to\infty,$ in $\Om^-.$  

\vskip0,3cm\noindent
This proves that $||u_k||$ is bounded in $\acca (\Om^+\cup\Om^-),$ and since $\Om^0$ is negligeable, 
this concludes this part of the proof. Then $u_k$  admits a subsequence weakly converging in $L^{2^*}.$ 

\vskip0,3cm
According to
(\ref{ps-diseq1}) and taking $v=u$ as a test function yields
\be\label{ps-diseq3}
\INT \Psi'(\na u_k)(\na u-\na u_k)\,dx\geq \lambda \INT u_k(u-u_k)\,dx+\INT b(x)(u_k)^{p-1}(u-u_k)\,dx +o(1)
\ee
so as $k\to \infty$ the right hand-side terms go to zero, and we obtain
\be\label{ps-diseq4}
\liminf \INT \Psi'(\na u_k)(\na u-\na u_k)\,dx\geq 0.
\ee
On the other hand, by convexity
\be\label{ps-convessita}
\INT \Psi(\na u)\,dx\geq \INT \Psi(\na u_k)\,dx + \INT \Psi'(\na u_k)(\na u-\na u_k)\,dx
\ee
So by (\ref{ps-diseq4}) and (\ref{ps-convessita})

\be\label{ps-convessita1}
\limsup\INT \Psi(\na u_k)\,dx\leq \limsup\left(\INT \Psi(\na u)\,dx - \INT \Psi'(\na u_k)(\na u-\na u_k)\,dx\right)
\ee
\[
\leq \INT \Psi(\na u)\,dx-\liminf \INT \Psi'(\na u_k)(\na u-\na u_k)\,dx\leq
\INT \Psi(\na u)\,dx
\]
By lower semicontinuity and convexity
\be\label{ps-convessita3}
\liminf\INT \Psi(\na u_k)\,dx\geq \INT \Psi(\na u)\,dx
\ee
We can conclude that
$$
\INT \Psi(\na u_k)\,dx\to \INT \Psi(\na u)\,dx.
$$
By Theorem \ref{natmod} $u_k$ admits a subsequence strongly converging in $L^{2^*},$ which concludes the proof of PS condition and of
Theorem \ref{esistenza}.
\par\medskip\noindent
\cvd
\par\medskip\noindent

\textbf{Proof of Theorem \ref{teocambiosegno}}
\par\medskip\noindent
We are now concerned with the existence of (possibly
sign-changing) nontrivial solutions $u$ of~(\p).
Let $(\lambda_k)$ denote the sequence of the eigenvalues of $-\Delta$ with
homogeneous Dirichlet condition, repeated according to multiplicity.
\par\medskip\noindent
Since the case $0<\lambda<\lambda_1$ is already contained in
Theorem~\ref{esistenza}, we may assume that $\lambda\geq \lambda_1$.
Let $k\geq 1$ be such that $\lambda_k \leq \lambda < \lambda_{k+1},$
$e_1,\ldots,e_k$ are eigenfunctions of $-\Delta,$ as defined in the introduction.
Finally, let $E_-=\mathrm{span}\{e_1,...,e_k\}$ and
$E_+=E_-^\perp$.
\par
Consider the functional $J$ defined in (\ref{funzionale-esteso})
We aim to apply the version of the Linking Theorem for convex functional presented by Szulkin in \cite{sz}.
Since
\[
\frac{\int_{\Om}\Psi(\na u)\,dx}{
\int_{\Om} |\na u|^2\,dx} \to \frac{1}{2}
\qquad\text{as $u\to 0$ in $\acca$,}
\]
as in the case $\Psi(\xi)=\frac{1}{2}|\xi|^2$
treated in \cite{gmm}, we deduce that there exist
$\rho>0$ and $\alpha>0$ such that $J(u)\geq\alpha$
whenever $u\in E_+$ with $\|u\|=\rho$.
On the other hand, there exists $e \in\acca\setminus E_-$
such that
\begin{gather*}
\lim_{\sublim{\|u\|\to\infty}{u\in \erre e\oplus E_-}}
J(u) = -\infty\,, \\
\end{gather*}
Again, this is proved in \cite{gmm} when
$\Psi(\xi)=\frac{1}{2}|\xi|^2$, but by $(\Psi_2)$ the
assertion is true also in our case.
Finally, it is clear that $J(u)\leq 0$ for every
$u\in E_-$.
\par
By the Linking type theorem in \cite{sz} (Theorem 3.4), there exist a PS sequence $(u_k)$ in $\acca$ and
we can continue, up to minor changes, as in the proof of Theorem~\ref{esistenza} to prove that there exists a subsequence of $(u_k)$
strongly converging in $L^{2^*}.$ This concludes the proof of Theorem \ref{teocambiosegno}, since the nontriviality of the solution
comes directly from the characterization of the critical level of the solution.
\par\medskip\noindent
\cvd
\par\medskip\noindent
\textbf{Acknowledgment} The author thanks Prof. Marco Degiovanni for very helpful conversations.

%%%%%%%%%%%%%%%%%%%%%%%%%%%%%%%%%%%%%%%%%%%%%%%%%%%%%%%%%%%%%%%%%%%%%%%%%%%%%
%%%%%%%%%%%%%%%%%%%%%%%%%%%%%%%%%%%%%%%%%%%%%%%%%%%%%%%%%%%%%%%%%%%%%%%%%%%%%

\end{document}